\documentclass{article}

\usepackage[french,british]{babel}

\usepackage[T1]{fontenc}

\usepackage{amsmath}

\usepackage{amsfonts}
\usepackage{amsbsy}
\usepackage{amssymb}

\usepackage{amsthm}

\theoremstyle{plain}

\theoremstyle{definition}
\newtheorem{definition}{Definition}

\newtheorem{lemma}{Lemma}
\newtheorem{theorem}{Theorem}
\newtheorem{corollary}{Corollary}

\theoremstyle{Remark}

\theoremstyle{definition}

\begin{document}

\title{Filtering of Continuous Time Periodically Correlated Isotropic Random Fields}

\author{Iryna Golichenko\thanks
{Department of Mathematical Analysis
and Probability Theory, National Technical University of Ukraine,
Kyiv 03056, Ukraine},
Oleksandr Masyutka\thanks
{Department of Mathematics and Theoretical Radiophysics,
Taras Shevchenko National University of Kyiv, Kyiv 01601, Ukraine},
Mikhail Moklyachuk\thanks
{Department of Probability Theory, Statistics and Actuarial
Mathematics, Taras Shevchenko National University of Kyiv, Kyiv 01601, Ukraine, Moklyachuk@gmail.com
}
     }

\date{\today}

\maketitle

\renewcommand{\abstractname}{Abstract}
\begin{abstract}
The problem of optimal linear estimation of functionals
depending on the unknown values of a random field
$\zeta(t,x)$, which is mean-square continuous periodically correlated
 with respect to  time argument $t\in\mathbb R$ and  isotropic on the unit sphere ${S_n}$ with respect to spatial argument $x\in{S_n}$.
Estimates are based on
observations of the field $\zeta(t,x)+\theta(t,x)$ at points $(t,x):t\leq 0,x\in S_{n}$,
 where $\theta(t,x)$ is an
uncorrelated with $\zeta(t,x)$
random field, which is mean-square continuous periodically correlated
 with respect to time argument $t\in\mathbb R$ and isotropic on the sphere ${S_n}$ with respect to spatial argument $x\in{S_n}$.
Formulas for calculating the mean square errors and the spectral characteristics of the optimal linear
estimate of functionals are derived in the case of spectral certainty where the spectral densities of the fields are exactly known.
Formulas that determine the least
favourable spectral densities and the minimax (robust) spectral
characteristics are proposed in the case where the spectral densities are not exactly known while a class of admissible spectral densities is given.
\end{abstract}

\vspace{2ex}

\vspace{2ex}

\section{Introduction}

Cosmological Principle (first coined by Einstein): the Universe is,
in the large, homogeneous and isotropic (J.~G.~Bartlett \cite{Bartlett}). Last
decades indicate growing interest to the
spatio-temporal data measured on the surface of a sphere. These data
includes cosmic microwave background (CMB) anisotropies (J.~G.~Bartlett \cite{Bartlett}, W.~Hu and S.~Dodelson \cite{Hu}, N.~Kogo and N.~Komatsu \cite{Kogo}, T.~Okamoto and W.~Hu \cite{Okamoto},
 P.~Adshead and W.~Hu \cite{Adshead}), medical imaging (R.~Kakarala \cite{Kakarala}), global
and land-based temperature data (P.~D.~Jones \cite{Jones}, T.~Subba Rao and G.~Terdik\cite{SubbaRao2006}), gravitational and geomagnetic data, climate model (G.~R.~North and R.~F.~Cahalan \cite{North}). Some basic results and references on the theory of isotropic random fields on a sphere can be found in the
books by M.~I.~Yadrenko \cite{Yadrenko} and A.~M.~Yaglom  \cite{Yaglom:1987a, Yaglom:1987b}. For more recent applications and results see new books by C.~Gaetan and X.~Guyon \cite{Gaetan},
N.~Cressie and C.~K.~Wikle \cite{Cressie}, D.~Marinucci and G.~Peccati \cite{Marinucci} and several
papers covering a number of problems in general for spatial temporal isotropic
observations (T.~Subba Rao and G.~Terdik\cite{SubbaRao2012}, G.~Terdik \cite{Terdik2015}).

Periodically correlated processes and fields are not homogeneous but have numerous properties similar to properties of
stationary processes and homogeneous fields. They describe appropriate models of numerous physical
and man-made processes. A comprehensive list of the existing
references up to the year 2005 on periodically correlated processes
and their applications was proposed by E.~Serpedin, F.~Panduru, I.~Sari and G.~B.~Giannakis\cite{Serpedin}. See also reviews by J.~Antoni \cite{Antoni} and A.~Napolitano \cite{Napolitano}. For more
details see a survey paper by W. A. Gardner \cite{Gardner} and book by H. L. Hurd and A. Miamee \cite{Hurd}. Note, that in the literature periodically
correlated processes are named in multiple different ways such as
cyclostationary, periodically nonstationary or cyclic correlated
processes.

The mean square optimal estimation problems for periodically
correlated with respect to time isotropic on a sphere random fields
are natural generalization of the linear extrapolation,
interpolation and filtering problems for stationary stochastic
processes and homogeneous random fields.
Effective methods of solution of the linear extrapolation, interpolation and filtering problems for stationary stochastic processes and random fields were developed under the condition of certainty where spectral densities of processes
and fields are known exactly (see, for example, selected works of A.~N.~Kolmogorov \cite{Kolmogorov}, survey article by T. Kailath \cite{Kailath}, books by Yu.~A.~Rozanov\cite{Rozanov}, N.~Wiener \cite{Wiener}, A.~M.~Yaglom  \cite{Yaglom:1987a, Yaglom:1987b}, M.~I.~Yadrenko \cite{Yadrenko}, articles by M. P. Moklyachuk and M. I. Yadrenko \cite{MoklyachukYadrenko:1979} - \cite{MoklyachukYadrenko:1980}).

The classical approach to the problems of interpolation, extrapolation and filtering of stochastic processes and random fields is based on the assumption that the spectral densities of processes and fields are known. In practice, however, complete information about the spectral density is impossible in most cases. To overcome this complication one finds parametric or nonparametric estimates of the unknown spectral densities or selects these densities by other reasoning. Then applies the classical estimation method provided that the estimated or selected density is the true one. This procedure can result in a significant increasing of the value of error as K.~S.~Vastola and H.~V.~Poor \cite{Vastola} have demonstrated with the help of some examples. This is a reason to search estimates which are optimal for all densities from a certain class of admissible spectral densities. These estimates are called minimax since they minimize the maximal value of the error of estimates.
Such problems arise when considering problems of automatic control theory, coding and signal processing in radar and sonar, pattern recognition problems of speech signals and images.
A comprehensive survey of results up to the year 1985 in minimax (robust) methods of data processing can be found in the paper by S.~A.~ Kassam and H.~V.~ Poor \cite{KassamPoor}. J. Franke \cite{Franke}, J. Franke and H. V. Poor \cite{Franke_Poor} investigated the minimax extrapolation  and filtering problems for stationary sequences with the help of convex optimization methods. This approach makes it possible to find equations that determine the least favorable spectral densities for different classes of densities.
 The paper by Ulf ~Grenander \cite{Grenander} should be marked as the first one where the minimax
approach to extrapolation problem for the functionals from stationary processes was
developed. For more details see, for example,
survey articles  M. P. Moklyachuk \cite{Moklyachuk:2000}, \cite{Moklyachuk:2001}, \cite{Moklyachuk:2015}
books by M.~Moklyachuk \cite{Moklyachuk:2008}, M.~Moklyachuk and O.~Masytka \cite{Moklyachuk:2012}, I. I. Golichenko and M. P. Moklyachuk \cite{Golichenko}.
In papers by I. I. Dubovets'ka, O.Yu. Masyutka and M.P. Moklyachuk\cite{Dubovetska1}, I. I.~Dubovets'ka and M. P.~ Moklyachuk \cite{Dubovetska4} - \cite{Dubovetska8}
the minimax-robust estimation problems (extrapolation, interpolation and filtering) are investigated
for linear functionals which depend on unknown values
of periodically correlated stochastic processes. Methods of
solution the minimax-robust estimation problems for time-homogeneous isotropic random fields on a sphere were developed by M. P. Moklyachuk \cite{Moklyachuk:1994} - \cite{Moklyachuk:1996}.
In papers by I. I. Dubovets'ka, O.Yu. Masyutka and M.P. Moklyachuk \cite{Dubovetska9} - \cite{Dubovetska11} results of investigation of minimax-robust estimation problems for periodically correlated isotropic random fields are
proposed.

In this article we deal with the problem of mean square optimal linear
estimation of the functional
\[
A\zeta ={\int_{0}^{\infty}}{\int_{S_n}} \,\,a(t,x)\zeta
(-t,x)\,m_n(dx)dt
\]
which depends on unknown values of a periodically correlated
(cyclostationary with period $T$) with respect to time isotropic on
the unit sphere ${S_n}$ in Euclidean space ${\mathbb E}^n$ random field
$\zeta(t,x)$, $t\le 0$, $x\in{S_n}$. Estimates are based on
observations of the field $\zeta(t,x)+\theta(t,x)$ at points
$(t,x)$, $t\le0$, $x\in{S_n}$, where $\theta(t,x)$ is an
uncorrelated with $\zeta(t,x)$ periodically correlated with respect
to time isotropic on the sphere ${S_n}$ random field. Formulas are
derived for computing the value of the mean-square error and the
spectral characteristic of the optimal linear estimate of the
functional $A\zeta$ in the case of spectral certainty, where spectral densities of the
fields are known. Formulas are proposed that determine the least
favourable spectral densities and the minimax-robust spectral
characteristic of the optimal estimate of the functional $A\zeta$
for concrete classes of spectral densities under the condition that
spectral densities are not known exactly while classes $D =D_f \times D_g$ of
admissible spectral densities are given.

\section{Spectral properties of periodically correlated isotropic on
a sphere random fields}
Let $S_{n}$ be a unit sphere in the $n$-dimensional Euclidean space ${\mathbb E}^n$, let $m_{n}(dx)$ be the Lebesgue measure on  $S_{n}$, and let
\[S_{m}^{l}(x),\, l=1, ... , h(m,n);\, m=0, 1, ...\]
be the orthonormal spherical harmonics of degree $m$, where $h(m,n)$ is the number of orthonormal spherical harmonics (see books by A. Erdelyi et al. \cite{Erdelyi} and C. M\"uller \cite{Muller} for more details).

A mean-square continuous random field ${\zeta(t,x)}$, $t\in\mathbb R$, $x\in{S_n}$, ${\zeta(t,x)}\in H=L_2(\Omega, \cal F, \mathbb P)$, where $L_2(\Omega, \cal F, \mathbb P)$ denotes the Hilbert space of random variables $\zeta$ with zero first moment, ${\mathbb E}{\zeta}=0$, and finite second moment, ${\mathbb E}|{\zeta}|^2<\infty$, is called periodically correlated (cyclostationary with period $T$) with respect to time isotropic on the sphere $S_{n}$ if for all $t,s\in\mathbb R$ and $x,y\in{S_n}$ the following property holds true
\[
{\mathbb E}\left({\zeta(t+T,x)} \overline
{\zeta (s+T,y)}\right)=B\left(t,s,\cos\vartheta \right),
\]
where $\cos\vartheta=(x,y)$, $\vartheta$ is the angular distance between points $x,y\in{S_n}$.

The correlation function $B\left(t,s,\cos\vartheta \right)$ of the mean-square continuous random field $\zeta(t,x)$ is continuous. It can be represented in the form of the series
\[
B\left(t,s,\cos\vartheta \right)=\frac{1}{\omega_n}\sum_{m=0}^{\infty}h(m,n) \frac
{C_{m}^{(n-2)/2}(\cos\vartheta)} {C_{m}^{(n-2)/2}(1)}\,\,
B^{\zeta}_m(t,s),
\]
where $\omega_n=(2\pi)^{n/2}\Gamma(n/2)$, $C_m^l(z)$ are the Gegenbauer polynomials (see book by M.~I.~Yadrenko \cite{Yadrenko}).

It follows from the Karhunen theorem that the random field $\zeta(t,x)$ itself can be represented in the form of the mean square convergent series (see K. Karhunen \cite{Karhunen}, I. I. Gikhman and A. V. Skorokhod \cite{Gikhman})
\begin{equation} \label{zeta1}
 {\zeta(t,x)}=
{\sum_{m=0}^{\infty}} {\sum_{l=1}^{h(m,n)}}  S_m^l (x)\zeta_{m}^l
(t),
\end{equation}
where
\[ \zeta_{m}^l (t)= {\int_{S_n}}{\zeta(t,x)} S_m^l
(x)\,m_n(dx).
\]
In this representation
\[\zeta_{m}^l (t),\,\, l=1,\ldots,h(m,n);\,\, t\in\mathbb R, m=0,1,\dots
\]
are mutually uncorrelated periodically correlated stochastic processes with the correlation functions $B^{\zeta}_m(t,s)$:
\[
{\mathbb E}\left(\zeta_{m}^l (t+T)\overline{\zeta_{u}^v
(s+T)}\right)=\delta_m^u \delta_l^v\,\,B^{\zeta}_m(t,s),\]
\[l,v=1,\ldots,h(m,n);\,\, m,u=0,1,\dots;\,\,t,s\in\mathbb R,
\]
where $\delta_l^v$ are the Kroneker delta-functions.

Consider two  mutually uncorrelated periodically correlated random fields  $\zeta(t,x)$ and $\theta(t,x)$. We construct the following sequences of stochastic functions
\begin{equation} \label{zj}
\{\zeta_{m}^{l}(j, u)=\zeta_{m}^{l}(u+jT),u\in [0,T), j\in\mathbb Z\},
\end{equation}
\begin{equation} \label{tj}
\{\theta_{m}^{l}(j, u)=\theta_{m}^{l}(u+jT),u\in [0,T), j\in\mathbb Z\}
\end{equation}
which correspond to the random fields $\zeta(t,x)$ and $\theta(t,x)$.
The sequences (\ref{zj}) and (\ref{tj})  form the $L_2([0,T);H)$-valued stationary sequences $\{\zeta_{m}^{l}(j), j\in\mathbb
Z\}$ and $\{\theta_{m}^{l}(j),j\in\mathbb Z\}$, respectively, with the correlation functions
\[
R_{m}^{\zeta}(k,j)=
\int_0^T{\mathbb E} [\zeta_{m}^{l}(u+kT)\overline{\zeta_{m}^{l}(u+jT)}]du
=\int_0^T B_{m}^{\zeta}(u+(k-j)T,u)du =
 R_{m}^{\zeta}(k-j), \]
 \[
R_{m}^{\theta}(k,j)= \int_0^T{\mathbb E} [\theta_{m}^{l}(u+kT)\overline{\theta_{m}^{l}(u+jT)}]du=
\int_0^T B_{m}^{\theta}(u+(k-j)T,u)du =
 R_{m}^{\theta}(k-j).\]

To describe properties of the stationary sequences $\{\zeta_{m}^{l}(j), j\in\mathbb
Z\}$ and $\{\theta_{m}^{l}(j),j\in\mathbb Z\}$ we define in the space $L_2([0,T);\mathbb{R})$ the following orthonormal basis
\[
\{\widetilde{e}_k=\frac{1}{\sqrt{T}}e^{2\pi
i\{(-1)^k\left[\frac{k}{2}\right]\}u/T}, k=1,2,\dots\}, \; \langle \widetilde{e}_j,\widetilde{e}_k\rangle=\delta_k^j.\]

Making use of the introduced basis the stationary sequences $\{\zeta_{m}^{l}(j), j\in\mathbb Z\}$ and $\{\theta_{m}^{l}(j),j\in\mathbb Z\}$ can be represented as follows
\begin{equation} \label{zeta}
\zeta_{m}^{l}(j)= \sum_{k=1}^\infty \zeta_{mk}^{l}(j)\widetilde{e}_k,
\end{equation}
\[\zeta_{mk}^{l}(j)=\langle\zeta_{m}^{l}(j),\widetilde{e}_k\rangle = \frac{1}{\sqrt{T}}
\int_0^T \zeta_{m}^{l}(j,v)e^{-2\pi i\{(-1)^k\left[\frac{k}{2}\right]\}v/T}dv,\]
\begin{equation} \label{teta}
\theta_{m}^{l}(j)= \sum_{k=1}^\infty \theta_{mk}^{l}(j)\widetilde{e}_k,
\end{equation}
\[
\theta_{mk}^{l}(j)=\langle\theta_{m}^{l}(j),\widetilde{e}_k\rangle= \frac{1}{\sqrt{T}}
\int_0^T \theta_{m}^{l}(j,v)e^{-2\pi i\{(-1)^k\left[\frac{k}{2}\right]\}v/T}dv.
\]

Components  of the constructed vector-valued stationary sequences $\{\zeta_{m}^{l}(j)=(\zeta_{mk}^{l}(j),k=1,2,\dots), j\in\mathbb Z\}$ and $\{\theta_{m}^{l}(j)=(\theta_{mk}^{l}(j),k=1,2,\dots),j\in\mathbb Z\}$ have the following properties \cite{Kallianpur},
\cite{Moklyachuk:1981}
\[
\mathbb{E}{\zeta_{mk}^{l}(j)}=0, \quad \|\zeta_{m}^{l}(j)\|^2_H=\sum_{k=1}^\infty
\mathbb{E}|\zeta_{mk}^{l}(j)|^2=R_{m}^{\zeta}(0), \quad
\mathbb{E}{\zeta_{mk}^{l}(j_1)}\overline{\zeta_{mn}^{l}(j_2)}=\langle
K_m^{\zeta}(j_1-j_2)e_k,e_n\rangle,
\]
\[
\mathbb{E}{\theta_{mk}^{l}(j)}=0, \quad \|\theta_{m}^{l}(j)\|^2_H=\sum_{k=1}^\infty
\mathbb{E}|\theta_{mk}^{l}(j)|^2=R_{m}^{\theta}(0), \quad
\mathbb{E}{\theta_{mk}^{l}(j_1)}\overline{\theta_{mn}^{l}(j_2)}=\langle
K_m^{\zeta}(j_1-j_2)e_k,e_n\rangle,
\]
where $\{e_k, k=1,2,\dots\}$ is a basis in the space $\ell_2$.
The correlation functions $K_m^{\zeta}(j)$  and $K_m^{\theta}(j)$ of the stationary sequences $\{\zeta_{m}^{l}(j)=(\zeta_{mk}^{l}(j),k=1,2,\dots), j\in\mathbb Z\}$ and $\{\theta_{m}^{l}(j)=(\theta_{mk}^{l}(j),k=1,2,\dots),j\in\mathbb Z\}$ are correlation operator functions in $\ell_2$.

The vector-valued stationary sequences $\{\zeta_{m}^{l}(j)=(\zeta_{mk}^{l}(j),k=1,2,\dots), j\in\mathbb Z\}$ and $\{\theta_{m}^{l}(j)=(\theta_{mk}^{l}(j),k=1,2,\dots),j\in\mathbb Z\}$ have the spectral density functions
\[F_{m}(\lambda)=\left\{ f_{m}^{kn}(\lambda)\right\}_{k,n = 1 }^{\infty}, \quad G_{m}(\lambda)=\left\{ g_{m}^{kn}(\lambda)\right\}_{k,n = 1 }^{\infty},\]
that are  operator-valued functions of variable $\lambda\in [-\pi,\pi)$ in the space $\ell_2$
if their correlation functions $K_m^{\zeta}(j)$  and $K_m^{\theta}(j)$ can be represented in the form
\[
\langle
K_m^{\zeta}(j)e_k,e_n\rangle
=\frac{1}{2\pi} \int _{-\pi}^{\pi} e^{ij\lambda}\langle F_{m}(\lambda) {e}_k,{e}_n\rangle d\lambda,
\]
\[
\langle
K_m^{\theta}(j)e_k,e_n\rangle
=\frac{1}{2\pi} \int _{-\pi}^{\pi} e^{ij\lambda}\langle G_{m}(\lambda) {e}_k,{e}_n\rangle d\lambda,
\]

For almost all   $\lambda\in [-\pi,\pi)$ the spectral densities $F_{m}(\lambda)$ and $G_{m}(\lambda)$  are kernel operators with integrable kernel norm

\[
\sum_{k=1}^\infty \frac{1}{2\pi} \int _{-\pi}^\pi \langle F_{m}(\lambda)
e_k, e_k\rangle d\lambda=\sum_{k=1}^\infty\langle K_m^{\zeta}(0)
e_k,e_k\rangle=\|\zeta_{m}^{l}(j)\|^2_H=R_{m}^{\zeta}(0),
\]
\[
\sum_{k=1}^\infty \frac{1}{2\pi} \int _{-\pi}^\pi \langle G_{m}(\lambda)
e_k, e_k\rangle d\lambda=\sum_{k=1}^\infty\langle K_m^{\theta}(0)
e_k,e_k\rangle=\|\theta_{m}^{l}(j)\|^2_H=R_{m}^{\theta}(0).
\]

\section{Hilbert space projection method of filtering}

Consider the problem of the mean square optimal linear estimation of the functional
\[
A\zeta ={\int_{0}^{\infty}}{\int_{S_n}} \,\,a(t,x)\zeta
(-t,x)\,m_n(dx)dt
\]
which depends on unknown values of a periodically correlated
with respect to time isotropic on
the unit sphere ${S_n}$ in Euclidean space ${\mathbb E}^n$ random field
$\zeta(t,x)$, $t\le 0$, $x\in{S_n}$.
Estimates are based on
observations of the field $\zeta(t,x)+\theta(t,x)$ at points
$(t,x)$, $t\le0$, $x\in{S_n}$, where $\theta(t,x)$ is an
uncorrelated with $\zeta(t,x)$ periodically correlated with respect
to time isotropic on the sphere ${S_n}$ random field.

It follows from representations \eqref{zeta1} that the functional $A\zeta$ can be represented in the form
\[A\zeta=\int_{0}^{\infty }\int_{S_{n}}a(t,x)\zeta(-t,x)m_{n}(dx)dt=\sum _{m=0}^{\infty }\sum _{l=1}^{h(m,n)}\int _{0}^{\infty }a_{m}^{l}(t)\zeta_{m}^{l} (-t)dt=\]
\[
=\sum _{m=0}^{\infty }\sum _{l=1}^{h(m,n)}\sum_{j=0}^\infty\int_{0}^{T} a_m^l(j,u)\zeta_m^l(-j,-u)du,
\]
\[ a_{m}^l (t)= {\int_{S_n}}{a(t,x)} S_m^l
(x)\,m_n(dx),
\]
\[
a_m^l(j,u)=a_m^l(u+jT),\,u\in [0,T),
\]
\[
\zeta_m^l(-j,-u)=\zeta_m^l(-u-jT),\, u\in [0,T).
\]

Taking into account the decomposition (\ref{zeta}) of stationary sequence $\{\zeta_m^l(j),j\in \mathbb Z\}$, the functional $A\zeta$ can be represented in the following form
\[A\zeta=\sum _{m=0}^{\infty }\sum _{l=1}^{h(m,n)}\sum_{j=0}^\infty\sum_{k=1}^{\infty} a_{mk}^l(j)\zeta_{mk}^l(-j)=
\sum _{m=0}^{\infty }\sum _{l=1}^{h(m,n)}\sum_{j=0}^\infty\vec{a}_m^l(j) ^{\top}\vec{\zeta}_m^l(-j),
\]
\[
\vec{\zeta}_m^l(-j)=(\zeta_{mk}^l(-j),k=1,2,\dots)^{\top},\]
\[
\vec{a}_m^l(j)=(a_{mk}^l(j),k=1,2,\dots)^{\top}=\]
\[
=(a_{m1}^l(j),a_{m3}^l(j),a_{m2}^l(j),\dots,a_{m(2k+1)}^l(j),a_{m(2k)}^l(j),\dots)^{\top},
\]
\[
 a_{mk}^l(j)=\langle a_m^l(j),\widetilde{e}_k\rangle =\frac{1}{\sqrt{T}} \int_0^T a_m^l(j,v)e^{-2\pi i\{(-1)^k\left[\frac{k}{2}\right]\}v/T}dv.
\]

\noindent We will assume that coefficients $\{\vec a_m^l(j),
j=0,1,\dots\}$ which form this representation satisfy the following conditions
\begin{equation} \label{aa}
\sum _{m=0}^{\infty }\sum _{l=1}^{h(m,n)}\sum_{j=0}^\infty \|
\vec a_m^l(j)\|<\infty,\quad \sum _{m=0}^{\infty }\sum _{l=1}^{h(m,n)}\sum_{j=0}^\infty (j+1)\|
\vec a_m^l(j)\|^2<\infty,\end{equation}
$$\begin{array} {c}\|\vec a_m^l(j)\|^2=\sum_{k=1}^\infty |a_{mk}^l(j)|^2.\end{array} $$

\noindent Under these condition the functional $A\zeta$ has finite second moment and operators defined below with the help of the coefficients $\{\vec a_m^l(j),
j=0,1,\dots\}$ are compact.

Denote by $L_{2}(F)$ the Hilbert space of complex vector functions
\[h(\lambda)=\left\{ h_{m}^l(\lambda):m=0,1, \ldots; l=1,2,\ldots , h(m,n) \right\},\; h_m^l(\lambda)=\left\{ h_{mk}^l \right\}_{k=1}^{\infty},\]
that satisfy condition
\[ \sum _{m=0}^{\infty }\sum _{l=1}^{h(m,n)}\int_{-\pi}^{\pi}(h_m^l(\lambda))^{\top}F_m(\lambda)\overline{h_m^l(\lambda)}d\lambda<\infty.\]

We denote by $L_{2}^{-}(F)$ the subspace of $L_{2}(F)$ generated by the functions
\[e^{ij\lambda}\delta_{k},\delta_{k}=\left\{ \delta_k^n \right\}_{n
= 1 }^{\infty}, k=1,2,\dots,  j\leq0,\]
where $\delta_k^k=1,\delta_k^n=0$, $k\neq n$.

Every linear estimate  $\hat{A}\zeta$ of the functional $A\zeta$ from observations of the sequence  $\{\zeta_m^l(j)+\theta_m^l(j), j\in \mathbb Z\}$ at points $j\leq0$  is defined by the spectral characteristic $h(\lambda)\in L_{2}^{-}(F+G)$ and is of the form
\begin{equation}  \label{est}
\hat{A}{\zeta}=\sum _{m=0}^{\infty }\sum _{l=1}^{h(m,n)}\int_{-\pi}^{\pi}(h_m^l(\lambda))^{\top}Z_m^{l\;\zeta+\theta}(d\lambda),
\end{equation}
where $Z_m^{l\;\zeta+\theta}(\Delta)=\{
Z_{mk}^{l\;\zeta+\theta}(\Delta) \}_{k = 1}^{\infty}$ is the orthogonal
stochastic measure of sum of sequences $\zeta_m^l(j)$ and $\theta_m^l(j)$.

Suppose that spectral densities of stationary sequence $\{\zeta_m^l(j)\}, \{\theta_m^l(j)\}$ admit the canonical factorizations (G.~Kallianpur and V.~Mandrekar \cite{Kallianpur}, M.~P.~Moklyachuk \cite{Moklyachuk:1981})
\begin{equation} \label{f}
F_m(\lambda)=\varphi_m(\lambda)(\varphi_m(\lambda))^*,\;
\varphi_m(\lambda)=\sum_{u=0}^\infty \varphi_m(u)e^{-iu\lambda},
\end{equation}
\begin{equation} \label{g}
G_m(\lambda)=\psi_m(\lambda)(\psi_m(\lambda))^*,\;
\psi_m(\lambda)=\sum_{u=0}^\infty \psi_m(u)e^{-iu\lambda},
\end{equation}
\begin{equation} \label{fg}
F_m(\lambda)+G_m(\lambda)=d_m(\lambda)(d_m(\lambda))^*,\; d_m(\lambda)=\sum_{u=0}^\infty d_m(u) e^{-iu\lambda},
\end{equation}
where matrices
\[d_m(u)=\{d_{mk}^r(u)\}^{r=\overline{1,M}}_{k=\overline{1,\infty}}, \quad\varphi_m(u)=\left\{\varphi_{mk}^r(u)
\right\}_{k=\overline{1,\infty}}^{r =\overline{1,M_1}}, \quad\psi_m(u)=\left\{ \psi_{mk}^r(u)
\right\}_{k=\overline{1,\infty}}^{r=\overline{1,M_2}}\]
 are coefficients of the canonical factorizations, $M_1$ is the multiplicity of $\zeta_m^l(j)$, $M_2$ is the multiplicity of $\theta_m^l(j)$ and $M$ is the multiplicity of $\zeta_m^l(j)+\theta_m^l(j)$.

The mean square error $\Delta(h;F,G)$ of  the linear estimate $\hat{A}\zeta$  with the spectral characteristic  $h_m^l(\lambda)=\sum_{j=0}^\infty \vec h_m^l(j) e^{-ij\lambda}$
can be represented in the form
\[\Delta(h;F,G)=E|A {\zeta}-\hat{A}
{\zeta}|^{2}=\]
\[= \sum _{m=0}^{\infty }\sum _{l=1}^{h(m,n)}(\|\mathbf{\Psi_m^l a_m^l}\|^2+\|\mathbf{D_m(a_m^l-h_m^l)}\|^2-\]
\[-\langle \mathbf{\Psi_m( a_m^l-h_m^l)},\mathbf{\Psi_m a_m^l}\rangle-\langle \mathbf{\Psi_m a_m^l},\mathbf{\Psi_m (a_m^l-h_m^l)}\rangle),\]
where operators $\mathbf{\Psi}$, $\mathbf{D}$ are defined as follows
\[\|\mathbf{\Psi_m a_m^l}\|^2=\sum_{q=0}^{\infty} \|(\mathbf{\Psi_m a_m^l})_{q}\|^2,\,
(\mathbf{\Psi_m a_m^l})_{q}=\sum_{j=0}^q(\psi_m(q-j))^{\top}\vec a_m^l(j),
\]
\[
\|\mathbf{D_m(a_m^l-h_m^l)}\|^2=\sum_{q=0}^{\infty}\|(\mathbf{D_m(a_m^l-h_m^l)})_{q}\|^2,
\]
\[
(\mathbf{D_m(a_m^l-h_m^l)})_{q}=\sum_{j=0}^q (d_m(q-j))^{\top}(\vec a_m^l(j)- \vec
h_m^l(j)),
\]
\[\langle \mathbf{\Psi_m(a_m^l-h_m^l)},\mathbf{\Psi_m a_m^l}\rangle=\sum_{q=0}^\infty\langle
(\mathbf{\Psi_m(a_m^l-h_m^l)})_{q},(\mathbf{\Psi_m a_m^l})_{q}\rangle.
\]

The spectral characteristic $h(F,G)$ of the optimal linear estimate $\hat A\zeta$ of the functional minimizes the value of the mean square error
\begin{equation} \label{zadacha}
\begin{split}
\Delta(F,G)=\Delta(h&(F,G);F,G)=\\
=\mathop {\min }\limits_{h \in
L_{2}^{-}(F+G)} \Delta (h;F,G)&=
\mathop {\min }\limits_{\hat{A}
{\zeta}}E|A{\zeta}-\hat{A} {\zeta}|^{2}.
\end{split}
\end{equation}

In the case where the spectral densities $G_m(\lambda)$ and $F_m(\lambda)+G_m(\lambda)$ admit factorizations
(\ref{g}) and (\ref{fg}),  the spectral characteristic $h(F,G)$, which is a solution of the optimization problem (\ref{zadacha}), and the mean square error $\Delta(F,G)$ of the optimal estimate $\hat A \zeta$  are determined by formulas
\begin{equation} \label{sg}
h_m^l(F,G)=A_m^l(\lambda)-(b_m(\lambda))^{\top}C_m^l(G)(\lambda),
\end{equation}
\begin{equation} \label{eg}
\Delta(F,G)=\sum _{m=0}^{\infty }\sum _{l=1}^{h(m,n)}\left[\|\mathbf{\Psi_m a_m^l}\|^2-\|\mathbf{B_m^* \Psi_m^*\Psi_m a_m^l}\|^2\right],
\end{equation}
where \[b_m(\lambda)=\{b_{mr}^k(\lambda)\}_{r=\overline{1,M}}^{k=\overline{1,\infty}},\quad
b_m(\lambda)=\sum_{u=0}^\infty b_m(u)e^{-iu\lambda},\quad
b_m(\lambda)d_m(\lambda)=I_{M},\]
\[C_m^l(G)(\lambda)=\sum_{j=0}^\infty (C_m^l(G))_{j}e^{-ij\lambda},\quad A_m^l(\lambda)=\sum_{j=0}^\infty \vec a_m^l(j) e^{-ij\lambda},\]
 \[(C_m^l(G))_{j}=(\mathbf{B^* \Psi^*\Psi a})_j=\sum_{q=0}^\infty\overline{b_m(q)}(\mathbf{\Psi_m^*\Psi_m a_m^l })_{j+q},\]
\[(\mathbf{\Psi_m^*\Psi_m a_m^l })_q=\sum_{u=0}^\infty
\overline{\psi_m(u)}(\mathbf{\Psi_m a_m^l})_{u+q},\]
\[\|\mathbf{B_m^* \Psi_m^*\Psi_m a_m^l}\|^2=\sum_{q=0}^\infty \|(\mathbf{B_m^* \Psi_m^*\Psi_m a_m^l})_q\|^2.\]

\indent In the case where the spectral densities $F_m(\lambda)$ and $F_m(\lambda)+G_m(\lambda)$ admit factorizations (\ref{f}) and (\ref{fg}), the spectral characteristic $h(F,G)$ and the mean square error $\Delta(F,G)$ of the optimal estimate $\hat A \zeta$  are defined by formulas
\begin{equation} \label{sf}
h_m^l(F,G)=(b_m(\lambda))^{\top}C_m^l(F)(\lambda),
\end{equation}
\begin{equation} \label{ef}
\Delta(F,G)=\sum _{m=0}^{\infty }\sum _{l=1}^{h(m,n)}\left[\|\mathbf{\Phi_m a_m^l}\|^2-\|\mathbf{B_m^* \Phi_m^*\Phi_m a_m^l}\|^2\right],
\end{equation}
\[C_m^l(F)(\lambda)=\sum_{j=0}^\infty (C_m^l(F))_{j}e^{-ij\lambda},\]
\[(C_m^l(F))_{j}=(\mathbf{B_m^* \Phi_m^*\Phi_m a_m^l})_j=\sum_{q=0}^\infty\overline{b_m(q)}(\mathbf{\Phi_m^*\Phi_m a_m^l})_{j+q}, \]
\[(\mathbf{\Phi_m^*\Phi_m a_m^l})_q=\sum_{u=0}^\infty\overline{\varphi_m(u)}(\mathbf{\Phi_m a_m^l})_{u+q},\]
\[(\mathbf{\Phi_m a_m^l})_{q}=\sum_{j=0}^q(\varphi_m(q-j))^{\top} \vec a_m^l(j),\quad \|\mathbf{\Phi_m a_m^l}\|^2=\sum_{q=0}^{\infty} \|(\mathbf{\Phi_m a_m^l})_{q}\|^2,
\]
\[\|\mathbf{B_m^* \Phi_m^*\Phi_m a_m^l}\|^2=\sum_{q=0}^\infty \|(\mathbf{B_m^*\Phi_m^*\Phi_m a_m^l})_q\|^2.
\]

Let us summarize our results and present them in the form of a theorem.

\begin{theorem}
Let  $\{\zeta(t,x), t\in \mathbb{R}, x\in S_n\}$ and $\{\theta(t,x), t\in \mathbb{R}, x\in S_n\}$ be mutually uncorrelated random fields, which are periodically correlated
 with respect to time argument $t\in\mathbb R$ and isotropic on the unit sphere ${S_n}$ with respect to spatial argument $x\in{S_n}$. Let the stationary sequences  $\{\zeta_m^l(j),j\in\mathbb Z\}$ and $\{\theta_m^l(j),j\in\mathbb Z\}$  constructed with the help of relations (\ref{zj}), (\ref{tj}), respectively, have spectral densities $F_m(\lambda)$ and $G_m(\lambda)$ that admit the canonical factorizations  (\ref{f}), (\ref{fg}) (or (\ref{g}), (\ref{fg})).  Let coefficients $\{\vec a_m^l(j), j=0,1,\dots\}$
that determine the functional $A\zeta$ satisfy conditions (\ref{aa}). Then the spectral characteristic $h(F,G)$  and the mean square error $\Delta(F,G)$ of the optimal estimate of the functional $A \zeta$ from observations of the field $\zeta(t,x)+\theta(t,x)$ at points $(t,x)$, $t\le0$, $x\in{S_n}$  are given by formulas (\ref{sf}), (\ref{ef}) (or
(\ref{sg}), (\ref{eg})), respectively. The optimal estimate $\hat A\zeta$ of the functional $A\zeta$ is calculated by the formula (\ref{est}).
\end{theorem}

\section{Minimax-robust method of filtering}

Formulas (\ref{sg}) --  (\ref{ef}) for calculating the spectral characteristic and the mean square error of the optimal linear estimate of the functional ${A} \zeta$ can be applied in
the case where the spectral densities $F_m(\lambda)$ and $G_m(\lambda)$ of stationary sequences $\{\zeta_m^l(j),j\in\mathbb Z\}$ and $\{\theta_m^l(j),j\in\mathbb Z\}$ constructed by relations (\ref{zj}), (\ref{tj}), are known.
 If the spectral densities are not  exactly known while a set of admissible densities $D=D_{F}\times D_{G}$ is specified, then the minimax approach to estimation of functional of unknown values is reasonable. That is we find the estimate which minimizes the mean square error for all spectral densities from a given set $D=D_{F}\times D_{G}$  simultaneously.

\begin{definition}
For a given class of spectral densities $D=D_{F}\times D_{G}$ the spectral densities $F_m^{0}(\lambda)\in D_{F}$ and $G_m^{0}(\lambda)\in D_{G}$ are called the least favorable in  $D$ for the optimal estimate of functional $A\zeta$ if
\[
\Delta(F^0,G^0)=\Delta(h(F^0,G^0);F^0,G^0)={\max_{\substack{(F,G)\in D}}} \Delta (h(F,G);F,G).
\]
\end{definition}
\begin{definition}
For a given class of spectral densities $D=D_{F}\times D_{G}$ the spectral characteristic $h^0(\lambda)$ of the optimal linear estimate of the functional $A \zeta$ is called minimax-robust if
the following relations hold true
$$\begin{array}{c}
 h^0(\lambda)\in H_{D}=\mathop {\bigcap}\limits_{(F,G)\in D} L_{2}^{-}(F+G),\end{array}$$
$$
 \mathop{\min_{\substack {h\in H_{D}}}} {\max_ {\substack {(F,G)\in D}}}\Delta(h;F,G)
=\mathop {\max_{\substack{(F,G)\in D}}}\Delta(h^0;F,G).
$$
\end{definition}

Taking into account the introduced definitions and relations (\ref{f}) -- (\ref{ef}) we can verify that the following lemma holds true.

\begin{lemma}
Spectral densities $F_m^{0}(\lambda)\in D_{F}$ and $G_m^{0}(\lambda)\in D_{G}$ which admit the canonical factorizations
(\ref{f}) - (\ref{fg}) are the least favorable in the class  $D=D_{F}\times D_{G}$  for the optimal linear estimation of the functional $A\zeta$ if the  coefficients of factorizations define a solution of the constrained optimization problem
\[
\Delta(F,G)=\sum _{m=0}^{\infty }\sum _{l=1}^{h(m,n)}\left[\|\mathbf{\Phi_m a_m^l}\|^2-\|\mathbf{B_m^* \Phi_m^*\Phi_m a_m^l}\|^2\right]\rightarrow sup,
\]
\begin{equation} \label{zf}
F_m(\lambda)=\varphi_m(\lambda) (\varphi_m(\lambda))^*\in D_{F},
\end{equation}
\[
G_m(\lambda)= d_m(\lambda) (d_m(\lambda))^*-\varphi_m(\lambda) (\varphi_m(\lambda))^* \in D_G,
\]
or the constrained optimization problem
\[
\Delta(F,G)=\sum _{m=0}^{\infty }\sum _{l=1}^{h(m,n)}\left[\|\mathbf{\Psi_m a_m^l}\|^2-\|\mathbf{B_m^* \Psi_m^*\Psi_m a_m^l}\|^2\right]\rightarrow sup,
\]
\begin{equation} \label{zg}
G_m(\lambda)=\psi_m(\lambda) (\psi_m(\lambda))^*\in D_{G},
\end{equation}
\[
F_m(\lambda)=d_m(\lambda) (d_m(\lambda))^*-\psi_m(\lambda) (\psi_m(\lambda))^* \in D_F.
\]
\end{lemma}

\begin{lemma}
Let the spectral density $F_m(\lambda)$ be given and admits the  factorization (\ref{f}). Then the spectral density $G_m^0(\lambda)$ is the least favorable in $D_G$ for the optimal estimation of the functional $A\zeta$  if
\[
F_m(\lambda)+G_m^0(\lambda)=d_m^0(\lambda) (d_m^0(\lambda))^*,
\]
where $d_m^0(\lambda)=\sum_{u=0}^\infty d_m^0(u) e^{-iu\lambda}$ and coefficients $\{ d_m^0(u), u=0,1,\dots\}$ are determined by solution of the constrained optimization problem
\begin{equation} \label{extrf}
\begin{split}
\sum _{m=0}^{\infty }\sum _{l=1}^{h(m,n)}\left[\|\mathbf{B_m^*\Phi_m^*\Phi_m  a_m^l}\|^2\right]&\rightarrow inf,\\
G_m(\lambda)=d_m(\lambda)  (d_m(\lambda))^*&-F_m(\lambda) \in D_G.
\end{split}
\end{equation}
\end{lemma}

\begin{lemma} \label{4.3}
Let the spectral density $G_m(\lambda)$  be given and admits the  factorization (\ref{g}). Then the spectral density
$F_m^0(\lambda)$ is the least favorable in $D_F$ for optimal estimation of the functional $A\zeta$ and admits canonical factorizations (\ref{f}), (\ref{fg}) if
\[
F_m^0(\lambda)+G_m(\lambda)=d_m^0(\lambda) (d_m^0(\lambda))^*,
\]
where $d_m^0(\lambda)=\sum_{u=0}^\infty d_m^0(u) e^{-iu\lambda}$ and coefficients $\{ d_m^0(u), u=0,1,\dots\}$ are determined by solution of the constrained optimization problem
\begin{equation} \label{extrg}
\begin{split}
\sum _{m=0}^{\infty }\sum _{l=1}^{h(m,n)}\left[\|\mathbf{B_m^*\Psi_m^*\Psi_m a_m^l}\|^2\right]&\rightarrow inf,\\
  F_m(\lambda)=d_m(\lambda)  (d_m(\lambda))^*&-G_m(\lambda) \in D_F.
\end{split}
\end{equation}
\end{lemma}

For more detailed analysis of properties of the least favorable spectral densities and the minimax-robust spectral characteristics we observe that the
 least favorable spectral densities $F_m^{0}(\lambda)\in D_{F}$, $G_m^{0}(\lambda)\in D_{G}$ and the minimax spectral characteristic $h^0=h(F^0,G^0)$ form a saddle point of the function  $\Delta(h;F,G)$ on the set $H_{D}\times D$. The saddle point inequalities
\[
\Delta(h^0;F,G)\leq\Delta(h^0;F^0,G^0)\leq\Delta(h;F^0,G^0), \]
\[
\forall
h\in H_{D}, \quad  \forall F\in D_{F}, \quad  \forall g\in D_{G}
\]
hold if $h^0=h(F^0,G^0)$, $h(F^0,G^0)\in H_{D}$ and $(F^0,G^0)$ is a solution of the constrained optimization problem
\begin{equation} \label{extr}
\Delta(h(F^0,G^0);F,G) \rightarrow {sup},  \quad (F,G)\in D,
\end{equation}
where the functional
\[
\Delta(h(F^0,G^0);F,G)=
\]
\[
=\sum _{m=0}^{\infty }\sum _{l=1}^{h(m,n)}\bigg[\frac{1}{2\pi}\int_{-\pi}^{\pi}
(C_m^l(G^0)(\lambda) ^{\top} b_m^0(\lambda)F_m(\lambda) \left ( b_m^0(\lambda) \right)^* \overline{C_m^l(G^0)(\lambda)}d\lambda+
\]
\begin{equation}\label{Delta}
+\frac{1}{2\pi}\int_{-\pi}^{\pi}(C_m^l(F^0)(\lambda)) ^{\top} b_m^0(\lambda)G_m(\lambda) \left ( b_m^0(\lambda) \right)^* \overline{C_m^l(F^0)(\lambda)}d\lambda\bigg].
\end{equation}

The constrained optimization problem (\ref{extr}) is equivalent to the following unconstrained optimization problem
\begin{equation}\label{eq12}
\Delta_D(F,G)=-\Delta(h(F^0,G^0);F,G)+\delta((F,G)|D)\rightarrow
inf,
\end{equation}
where $\delta((F,G)|D)$ is the indicator function of the set $D$. Solution $(F^{0}(\lambda), G^0(\lambda))$ to the extremum problem (\ref{eq12}) is determined by the condition $0\in\partial \Delta_D(F^{0}, G^{0})$  which is necessary for the point  $(F^{0},G^0)$ to belong to the  set of minimums  of a convex functional. Here
$\partial \Delta_D(F^{0}, G^{0})$ is a subdifferential of the convex functional $\Delta_D(F,G)$ at point $(F,G)=(F^{0}, G^{0})$ (see  R. T. Rockafellar \cite{Rockafellar}, M.~P.~Moklyachuk \cite{Moklyachuk:2008nonsm}).

The form (\ref{Delta}) of the functional $\Delta (h(F^{0} ,G^{0} );F,G)$  is convenient for application the method of Lagrange multipliers for
finding solution to the problem (\ref{eq12}).
Making use the method of Lagrange multipliers and the form of
subdifferentials of the indicator functions $\delta((F,G)|D)$
we describe relations that determine the least favourable spectral densities in some special classes
of spectral densities (see books by M.~Moklyachuk \cite{Moklyachuk:2008}, M.~Moklyachuk and O.~Masytka \cite{Moklyachuk:2012}, I. I. Golichenko and M. P. Moklyachuk \cite{Golichenko} for more details).

\section{The least favorable spectral densities in the class $D_{0}\times D_V^U$}

Consider the problem of minimax estimation of the functional $A\zeta$
depending on the unknown values of the random field
$\{\zeta(t,x), t\in \mathbb R, x\in S_n\}$, which is periodically correlated
 with respect to the time argument $t\in\mathbb R$ and isotropic on the sphere ${S_n}$ with respect to spatial argument $x\in{S_n}$
 based on observations of the random field
$\zeta(t,x)+\theta(t,x)$ at points $(t,x):$ $ t\leq0, x\in{S_n},$
under the condition that  spectral densities
 $F_m(\lambda)$, $G_m(\lambda)$ of stationary sequences $\{ \zeta_m^l(j), j\in \mathbb{Z}\}$ and $\{ \theta_m^l(j), j\in \mathbb{Z}\}$ which are constructed with the help of relations (\ref{zj}), (\ref{tj}), respectively,
 are not known exactly while there are specified the
following pairs of sets of admissible spectral densities.

The first pair is

\[D_{0}^{1} =\left\{F(\lambda )|\, \frac{1}{2\pi\omega_{n} }\sum_{m=0}^{\infty}h(m,n) \int _{-\pi }^{\pi }\text{Tr}\,  F_{m}(\lambda )d\lambda =p \right\},\]
\[ {D_{V}^{U}} ^{1} =\bigg\{G(\lambda )|Tr V_{m} (\lambda )\le \text{Tr}\, G_{m}(\lambda )\le Tr U_{m} (\lambda ), \] \[\frac{1}{2\pi\omega_{n} }\sum_{m=0}^{\infty}h(m,n) \int _{-\pi }^{\pi } \text{Tr}\,  G_{m}(\lambda )d\lambda  =q\bigg\}.\]

The second pair of sets of admissible spectral densities is
\[D_{0}^{2} =\biggl\{F(\lambda )|\frac{1}{2\pi \omega_{n}}\sum_{m=0}^{\infty}h(m,n) \int _{-\pi }^{\pi }F_{m}^{kk} (\lambda )d\lambda =p_{k}, k=1,2,\dots \biggr\},\]
\[{D_{V}^{U}} ^{2} =\biggl\{G(\lambda )|V_{m}^{kk} (\lambda )\le G_{m}^{kk} (\lambda )\le U_{m}^{kk} (\lambda ), \]
\[ \frac{1}{2\pi\omega_{n} }\sum_{m=0}^{\infty}h(m,n) \int _{-\pi }^{\pi }G_{m}^{kk} (\lambda )d\lambda  =q_{k}, k=1,2,\dots\biggr\}.\]

The third pair of sets of admissible spectral densities is
\[D_{0}^{3} =\left\{F(\lambda )|\frac{1}{2\pi \omega_{n}}\sum_{m=0}^{\infty}h(m,n) \int _{-\pi }^{\pi }\left\langle B, F_{m}(\lambda )\right\rangle d\lambda  =p\right\},\]
\[{D_{V}^{U}} ^{3} =\biggl\{G(\lambda )|\left\langle B_{2} , V_{m} (\lambda )\right\rangle \le \left\langle B_{2} , G_{m}(\lambda )\right\rangle \le \left\langle B_{2} , U_{m} (\lambda )\right\rangle , \] \[\frac{1}{2\pi\omega_{n} }\sum_{m=0}^{\infty}h(m,n) \int _{-\pi }^{\pi }\left\langle B_{2} , G_{m}(\lambda )\right\rangle d\lambda  =q\biggr\}.\]

The forth pair of sets of admissible spectral densities is
\[D_{0}^{4} =\left\{F(\lambda )|\, \frac{1}{2\pi \omega_{n}}\sum_{m=0}^{\infty}h(m,n) \int _{-\pi }^{\pi }F_{m}(\lambda )d\lambda =P \right\},\]
\[{D_{V}^{U}} ^{4} =\biggl\{G(\lambda )|V_{m} (\lambda )\le G_{m}(\lambda )\le U_{m} (\lambda ), \]
\[\frac{1}{2\pi \omega_{n}}\sum_{m=0}^{\infty}h(m,n) \int _{-\pi }^{\pi }G_{m}(\lambda )d\lambda  =Q\biggr\}.\]

Here $V_{m} (\lambda ), U_{m} (\lambda )$ are given matrices of spectral densities, $p, q, p_{k}, q_{k}, k=1,2,\dots$ are given numbers, $B_{1}, B_{2}, P, Q$ are given positive-definite Hermitian matrices.

From the condition $0\in \partial \Delta _{D} (F^{0} ,G^{0} )$ we find the following equations which determine the least favourable spectral densities for these given sets of admissible spectral densities.

For the first pair $D_{0}^{1}\times {D_{V}^{U}} ^{1}$ we have equations
\begin{equation} \label{GrindEQ__1_20_}
\sum_{l=1}^{h(m,n)}C_{m}^{l}(G^{0})(\lambda)(C_{m}^{l}(G^{0})(\lambda))^{*}=\alpha_{m} ^{2} d_{m}^{0} (\lambda )^{\top} \overline{d_{m} ^{0} (\lambda )},
\end{equation}
\begin{equation} \label{GrindEQ__1_21_}
\sum_{l=1}^{h(m,n)}C_{m}^{l}(F^{0})(\lambda)(C_{m}^{l}(F^{0})(\lambda))^{*}=(\beta_{m}^{2} +\gamma _{m_{1}} (\lambda )+\gamma _{m_{2}} (\lambda ))d_{m}^{0} (\lambda )^{\top} \overline{d_{m}^{0} (\lambda )},
\end{equation}
where
\[\gamma _{m_{1}} (\lambda )\le 0\,\, \text{and}\,\, \gamma _{m_{1}} (\lambda )=0\,\,\text{if}\,\, \text{Tr}\,  G_{m}^{0} (\lambda )>\text{Tr}\,  V_{m}(\lambda ),\]
\[\gamma _{m_{2}} (\lambda )\ge 0\,\, \text{and}\,\, \gamma _{m_{2}} (\lambda )=0\,\,\text{if}\,\, \text{Tr}\,  G_{m}^{0} (\lambda )<\text{Tr}\,  U_{m}(\lambda ),\]
and $\alpha_{m} ^{2} , \beta_{m} ^{2}$ are unknown Lagrange multipliers.

For the second pair $D_{0}^{2}\times {D_{V}^{U}} ^{2}$ we have equations
\begin{equation} \label{GrindEQ__1_22_}
\sum_{l=1}^{h(m,n)}C_{m}^{l}(G^{0})(\lambda)(C_{m}^{l}(G^{0})(\lambda))^{*}=d_{m}^{0} (\lambda )^{\top} \left\{\alpha _{mk}^{2} \delta _{kl} \right\}_{k,l=1}^{\infty} \overline{d_{m}^{0} (\lambda )},
\end{equation}
\begin{equation} \label{GrindEQ__1_23_}
\sum_{l=1}^{h(m,n)}C_{m}^{l}(F^{0})(\lambda)(C_{m}^{l}(F^{0})(\lambda))^{*}=d_{m}^{0} (\lambda )^{\top} \left\{(\beta _{mk}^{2} +\gamma _{m_{1}k} (\lambda )+\gamma _{m_{2}k} (\lambda ))\delta _k^l \right\}_{k,l=1}^{\infty} \overline{d_{m}^{0} (\lambda )},
\end{equation}
where
\[ \gamma _{m_{1}k} (\lambda )\le 0\,\,\text{and}\,\, \gamma _{m_{1}k} (\lambda )=0\,\,\text{if}\,\, G^{0kk}_{m} (\lambda )>V^{kk}_{m}(\lambda ),\]
\[ \gamma _{m_{2}k} (\lambda )\ge 0\,\,\text{and}\,\,\gamma _{m_{2}k} (\lambda )=0\,\,\text{if}\,\, G^{0kk}_{m} (\lambda )<U^{kk}_{m}(\lambda ),\]
and $\alpha _{mk}^{2} , \beta _{mk}^{2}$ are unknown Lagrange multipliers.

For the third pair $D_{0}^{3}\times {D_{V}^{U}} ^{3}$ we have equations
\begin{equation} \label{GrindEQ__1_24_}
\sum_{l=1}^{h(m,n)}C_{m}^{l}(G^{0})(\lambda)(C_{m}^{l}(G^{0})(\lambda))^{*}=\alpha_{m} ^{2} d_{m}^{0} (\lambda )^{\top}B_{1} \overline{d_{m}^{0} (\lambda )},
\end{equation}
\begin{equation} \label{GrindEQ__1_25_}
\sum_{l=1}^{h(m,n)}C_{m}^{l}(F^{0})(\lambda)(C_{m}^{l}(F^{0})(\lambda))^{*}=(\beta_{m}^{2} +\gamma _{m_{1}} (\lambda )+\gamma _{m_{2}} (\lambda ))d_{m}^{0} (\lambda )^{\top} B_{2} \overline{d_{m}^{0} (\lambda )};
\end{equation}
where
\[\gamma _{m_{1}} (\lambda )\le 0\,\,\text{and}\,\, \gamma _{m_{1}} (\lambda )=0\,\,\text{if}\,\, \left\langle B_{2} , G_{m}^{0} (\lambda )\right\rangle >\left\langle B_{2} , V_{m} (\lambda )\right\rangle ,\]
\[ \gamma _{m_{2}} (\lambda )\ge 0\,\,\text{and}\,\,\gamma _{m_{2}} (\lambda )=0 \,\,\text{if}\,\, \left\langle B_{2} , G_{m}^{0} (\lambda )\right\rangle <\left\langle B_{2} , U_{m} (\lambda )\right\rangle ,\]
and $\alpha_{m} ^{2} , \beta_{m} ^{2}$ are unknown Lagrange multipliers.

For the forth pair $D_{0}^{4}\times {D_{V}^{U}} ^{4}$ we have equations
\begin{equation} \label{GrindEQ__1_26_}
\sum_{l=1}^{h(m,n)}C_{m}^{l}(G^{0})(\lambda)(C_{m}^{l}(G^{0})(\lambda))^{*}=d_{m}^{0} (\lambda )^{\top} \vec{\alpha_{m}}\cdot \vec{\alpha_{m}}^{*}\overline{d_{m}^{0} (\lambda )},
\end{equation}
\begin{equation} \label{GrindEQ__1_27_}
\sum_{l=1}^{h(m,n)}C_{m}^{l}(F^{0})(\lambda)(C_{m}^{l}(F^{0})(\lambda))^{*}=
d_{m}^{0} (\lambda )^{\top} (\vec{\beta}\cdot \vec{\beta}^{*}+\Gamma _{m_{1}} (\lambda )+\Gamma _{m_{2}} (\lambda ))\overline{d_{m}^{0} (\lambda )}.
\end{equation}
where
$\Gamma _{m_{1}} (\lambda ), \Gamma _{m_{2}} (\lambda )$ are Hermitian matrices,
\[\Gamma _{m_{1}} (\lambda )\le 0\,\,\text{and}\,\, \Gamma _{m_{1}} (\lambda )=0\,\,\text{if}\,\, G_{m}^{0} (\lambda )>V_{m}(\lambda ),\]
\[\Gamma _{m_{2}} (\lambda )\ge 0\,\,\text{and}\,\, \Gamma _{m_{2}} (\lambda )=0 \,\,\text{if}\,\, G_{m}^{0} (\lambda )<U_{m} (\lambda ),\]
and $   \vec{\alpha_{m} }, \vec{\beta_{m} }$ are unknown Lagrange multipliers.

\begin{theorem}
The least favorable spectral densities  $F_m^0(\lambda)$, $G_m^0(\lambda)$  in the classes $D_0\times D_V^U$  for the optimal estimate of the functional $A\zeta$  are determined by relations
\eqref{GrindEQ__1_20_}, \eqref{GrindEQ__1_21_} for the first pair $D_{0}^{1}\times {D_{V}^{U}} ^{1}$ of sets of admissible spectral densities
(\eqref{GrindEQ__1_22_}, \eqref{GrindEQ__1_23_} for the second  pair $D_{0}^{2}\times {D_{V}^{U}} ^{2}$ of sets of admissible spectral densities,
\eqref{GrindEQ__1_24_}, \eqref{GrindEQ__1_25_} for the third  pair $D_{0}^{3}\times {D_{V}^{U}} ^{3}$ of sets of admissible spectral densities,
\eqref{GrindEQ__1_26_}, \eqref{GrindEQ__1_27_} for the fourth pair $D_{0}^{4}\times {D_{V}^{U}} ^{4}$ of sets of admissible spectral densities),
factorizations (\ref{f}), (\ref{g}),  (\ref{fg}), constrained optimization problem  (\ref{zf}) or (\ref{zg}), and restrictions  on densities from
the corresponding classes $D_0\times D_V^U$. The minimax spectral characteristic $h(F^0,G^0)$ of the optimal estimate $\hat A \zeta$ is calculated by (\ref{sf}) or (\ref{sg}). The mean square error $\Delta(F^0,G^0)$ is calculated by  (\ref{ef}) or (\ref{eg}).
\end{theorem}

In the case where one of spectral densities $F_m(\lambda)$ or $G_m(\lambda)$ from the corresponding classes is known we have the following corollary from the theorem.

\begin{corollary}
If the spectral density $F_m(\lambda)\in D_0$ is known and admits the canonical factorization  (\ref{f}), then the least favorable spectral densities $G_m^0(\lambda)$ in the classes ${D_V^U}^{k}$, $k=1,2,3,4$ are determined by relations (\ref{g}), (\ref{fg}),  (\ref{extrf}), equations \eqref{GrindEQ__1_21_}, \eqref{GrindEQ__1_23_}, \eqref{GrindEQ__1_25_}, \eqref{GrindEQ__1_27_}
correspondingly to $k=1,2,3,4$ and by  restrictions on densities from classes ${D_V^U}^{k}$, $k=1,2,3,4$.
If the spectral density $G_m(\lambda)\in D_V^U$ is known and admits the canonical factorization  (\ref{g}), then the least favorable spectral densities $F_m^0(\lambda)$ in the classes $D_0^k$, $k=1,2,3,4$ are determined by relations (\ref{f}), (\ref{fg}),  (\ref{extrg}), equations \eqref{GrindEQ__1_20_},  \eqref{GrindEQ__1_22_},  \eqref{GrindEQ__1_24_},  \eqref{GrindEQ__1_26_})
correspondingly to $k=1,2,3,4$ and by  restrictions on densities from classes ${D_V^U}^{k}$, $k=1,2,3,4$.
 The minimax spectral characteristic $h(F^0,G^0)$ of the optimal estimate $\hat A \zeta$ is calculated by (\ref{sf}) or (\ref{sg}). The mean square error $\Delta(F^0,G^0)$ is calculated by (\ref{ef}) or (\ref{eg}).
\end{corollary}

\section{Conclusions}
In this paper we propose formulas for calculating the mean square error and the spectral characteristic of the optimal linear estimate of the
functional
\[
A\zeta ={\int_{0}^{\infty}}{\int_{S_n}} \,\,a(t,x)\zeta
(-t,x)\,m_n(dx)dt
\]
depending on unknown values of a mean-square continuous periodically correlated
(cyclostationary with period $T$) with respect to time argument and isotropic on
the unit sphere ${S_n}$ in Euclidean space ${\mathbb E}^n$ random field
$\zeta(t,x)$, $t\in\mathbb R$, $x\in{S_n}$. Estimates are based on
observations of the field $\zeta(t,x)+\theta(t,x)$ at points
$(t,x)$, $t\le0$, $x\in{S_n}$, where $\theta(t,x)$ is an
uncorrelated with $\zeta(t,x)$ mean-square continuous periodically correlated with respect
to time argument and isotropic on the sphere ${S_n}$ random field.
The problem is investigated in the case of spectral certainty where matrices of spectral densities of random fields are known exactly and in the case of spectral uncertainty where matrices of spectral densities of random fields are not known exactly while some classes of admissible spectral density matrices are given.
 We derive formulas for calculation the spectral characteristic and the mean-square error of the optimal linear estimate of the functional $A\zeta$
in the case of spectral certainty, where spectral densities $F_{m}(\lambda ), G_{m}(\lambda )$ of the stationary sequences that generate the random fields $\zeta(t,x)$, $\theta (t,x)$ are known exactly.

 We propose a representation of the mean square
error in the form of a linear functional in the $L_1\times L_1$ space with
respect to spectral densities $(F,G)$, which allows us to solve the
corresponding constrained optimization problem and describe the minimax
(robust) estimates of the functional $A\zeta$
for concrete classes of spectral densities under the condition that
spectral densities are not known exactly while classes $D =D_f \times D_g$ of
admissible spectral densities are given.


\begin{thebibliography}{99}


\bibitem{Adshead}
\newblock P. Adshead, and W. Hu,
\newblock \emph{Fast computation of first-order feature-bispectrum corrections},
\newblock Phys. Rev., vol. D85, 103531, 2012.

\bibitem{Antoni}
\newblock J. Antoni,
\newblock \emph{Cyclostationarity by examples},
\newblock Mechanical Systems and Signal Processing, vol. 23, pp. 987--1036, 2009.

\bibitem{Bartlett}
\newblock J. G. Bartlett,
\newblock \emph{The standard cosmological model and cmb anisotropies},
\newblock New Astron. Rev., vol. 43, pp. 83--109, 1999.

\bibitem{Cressie}
\newblock N. Cressie, and C. K. Wikle,
\newblock \emph{Statistics for spatio-temporal data},
\newblock Wiley Series in Probability and Statistics, 2011.


\bibitem{Dubovetska1}
\newblock I. I. Dubovets'ka, O.Yu. Masyutka, and M.P. Moklyachuk,
\newblock \emph{Interpolation of periodically correlated stochastic sequences},
\newblock Theory of Probability and Mathematical Statistics, vol. 84, pp. 43--56, 2012.

\bibitem{Dubovetska4}
\newblock I. I. Dubovets'ka, and M. P. Moklyachuk,
\newblock \emph{Filtration of linear functionals of periodically correlated sequences},
\newblock Theory of Probability and Mathematical Statistics, vol. 86, pp. 51--64, 2013.

\bibitem{Dubovetska6}
\newblock I. I. Dubovets'ka, and M. P. Moklyachuk,
\newblock \emph{Extrapolation of periodically correlated processes from observations with noise},
\newblock Theory of Probability and Mathematical Statistics, vol. 88, pp. 43--55, 2013.

\bibitem{Dubovetska7}
\newblock I. I. Dubovets'ka, and M. P. Moklyachuk,
\newblock \emph{Minimax estimation problem for periodically correlated stochastic processes},
\newblock Journal of Mathematics and System Science, vol. 3, no. 1, pp. 26--30, 2013.

\bibitem{Dubovetska8}
\newblock I. I. Dubovets'ka, and M. P. Moklyachuk,
\newblock \emph{On minimax estimation problems for periodically correlated stochastic processes},
\newblock Contemporary Mathematics and Statistics, vol.2, no. 1, pp. 123--150, 2014.


\bibitem{Dubovetska9}
\newblock I. I. Dubovets'ka, O. Yu. Masyutka, and M. P. Moklyachuk,
\newblock \emph{Filtering problems for periodically correlated isotropic random fields},
\newblock Mathematics and Statistics, vol.2, no. 4, pp. 162--171, 2014.

\bibitem{Dubovetska10}
\newblock I. I. Dubovets'ka, O. Yu. Masyutka, and M. P. Moklyachuk,
\newblock \emph{Estimation problems for periodically correlated isotropic random fields},
\newblock Methodology and Computing in Applied Probability, vol.17, no. 1, pp. 41--57, 2015.


\bibitem{Dubovetska11}
\newblock I. I. Dubovets'ka, O. Yu. Masyutka, and M. P. Moklyachuk,
\newblock \emph{Minimax-robust fitering of functionals from periodically correlated random fields},
\newblock Cogent Mathematics, vol.2, 1074327, 2015.

\bibitem{Erdelyi}
\newblock A. Erdelyi, W. Magnus, F. Oberhettinger, F. G. Tricomi,
\newblock \emph{Higher transcendental functions. Vol. II.}
\newblock Bateman Manuscript Project. New York-Toronto-London: McGraw-Hill Book Co., Inc. XVII, 1953.



\bibitem{Franke}
\newblock J. Franke,
\newblock \emph{Minimax robust prediction of discrete time series},
\newblock Z. Wahrscheinlichkeitstheor. Verw. Gebiete, vol. 68, pp. 337--364, 1985.

\bibitem{Franke_Poor}
\newblock J. Franke and H. V. Poor,
\newblock \emph{Minimax-robust filtering and finite-length robust predictors},
\newblock Robust and Nonlinear Time Series Analysis. Lecture Notes in Statistics, Springer-Verlag,
vol. 26, pp. 87--126, 1984.


\bibitem{Gaetan}
\newblock C. Gaetan, and X. Guyon,
\newblock \emph{Spatial statistics and modeling},
\newblock Springer Series in Statistics, vol. 81, Springer Science+Business Media, 2010.

\bibitem{Gikhman}
 \newblock I. I. Gikhman and A. V. Skorokhod,
\newblock \emph{The theory of stochastic processes. I.},
\newblock Berlin: Springer, 2004.

\bibitem{Gardner}
\newblock W. A. Gardner,
\newblock \emph{Cyclostationarity in communications and signal processing},
\newblock  New York: IEEE Press, 1994.

\bibitem{Gladyshev}
\newblock E. G. Gladyshev,
\newblock \emph{Periodically correlated random sequences},
\newblock Sov. Math. Dokl. vol. 2, pp. 385--388, 1961.


\bibitem{Golichenko}
\newblock I. I. Golichenko and M. P. Moklyachuk,
\newblock \emph{Estimates of functionals of periodically correlated processes},
\newblock Kyiv: NVP ``Interservis", 2014.

\bibitem{Grenander}
\newblock U. Grenander,
\newblock \emph{A prediction problem in game theory},
\newblock Arkiv f\"or Matematik, vol. 3, pp. 371--379, 1957.


\bibitem{Hu}
\newblock W. Hu, and S. Dodelson,
\newblock \emph{Cosmic microwave background anisotropies},
\newblock Annual Review of Astronomy and Astrophysics, vol. 40, pp. 171--216, 2002.

\bibitem{Hurd}
\newblock H. L. Hurd, and A. Miamee,
\newblock \emph{Periodically Correlated random sequences: Spectral theory and practice},
\newblock Wiley Series in Probability and Statistics; Wiley Interscience. Hoboken, NJ: John Wiley and Sons, 2007.

\bibitem{Jones}
\newblock P. D. Jones,
\newblock \emph{Hemispheric surface air temperature variations: A reanalysis and an update to 1993},
\newblock Journal of Climate, vol. 7, pp. 1794-1802, 1994.

\bibitem{Kailath}
\newblock T. Kailath,
\newblock \emph{A view of three decades of linear filtering theory},
\newblock IEEE Transactions on Information Theory, Vol. 20, pp. 146--181, 1974.

\bibitem{Kakarala}
\newblock R. Kakarala,
\newblock \emph{The bispectrum as a source of phase-sensitive invariants for Fourier descriptors: A group-theoretic approach},
\newblock Journal of Mathematical Imaging and Vision, vol. 44, pp. 341--353, 2012.


\bibitem{Kallianpur}
\newblock G. Kallianpur, and V. Mandrekar,
\newblock \emph{Spectral theory of stationary H-valued processes},
\newblock J.  Multivariate Analysis, vol. 1, pp. 1--16, 1971.

\bibitem{Karhunen}
\newblock K. Karhunen,
\newblock \emph{Uber lineare Methoden in der Wahrscheinlichkeitsrechnung},
\newblock Annales Academiae Scientiarum Fennicae. Ser. A I, no. 37, 1947.

\bibitem{KassamPoor}
\newblock S. A. Kassam, and H. V. Poor,
\newblock \emph{Robust techniques for signal processing: A survey},
\newblock Proceedings of the IEEE, vol. 73, no. 3, pp. 433--481, 1985.

\bibitem{Kogo}
\newblock N. Kogo, and N.Komatsu,
\newblock \emph{Angular trispectrum of cmb temperature anisotropy from primordial non-Gaussianity with the full radiation
transfer function},
\newblock Phys. Rev., vol. D73, pp. 083007--083012, 2006.

\bibitem{Kolmogorov}
\newblock A. N. Kolmogorov,
\newblock \emph{Selected works by A. N. Kolmogorov. Vol. II: Probability theory and mathematical statistics. Ed. by A. N. Shiryayev},
\newblock Mathematics and its Applications. Soviet Series. 26. Dordrecht etc.
Kluwer Academic Publishers, 1992.

\bibitem{Marinucci}
\newblock D. Marinucci, and G. Peccati,
\newblock \emph{Random Fields on the sphere},
\newblock London Mathematical Society Lecture Notes Series, vol. 389, Cambridge
University Press, Cambridge, 2011.


\bibitem {Moklyachuk:1981}
\newblock M. P. Moklyachuk,
\newblock \emph{Estimation of linear functionals of stationary stochastic processes and a two-person zero-sum game},
\newblock Stanford University Technical Report, no. 169, 1981.


\bibitem{Moklyachuk:1994}
\newblock M. P. Moklyachuk,
\newblock \emph{Minimax filtering of time-homogeneous isotropic random fields on a sphere},
\newblock Theory of Probability and Mathematical Statistics, vol. 49, pp. 137--146, 1994.

\bibitem{Moklyachuk:1995}
\newblock M. P. Moklyachuk,
\newblock \emph{Extrapolation of time-homogeneous random fields that are isotropic on a sphere. I},
\newblock Theory of Probability and Mathematical Statistics, vol. 51, pp. 137--146 , 1995.


\bibitem{Moklyachuk:1996}
\newblock M. P. Moklyachuk,
\newblock \emph{Extrapolation of time-homogeneous random fields that are isotropic on a sphere. II},
\newblock Theory of Probability and Mathematical Statistics, vol.53, pp. 137--148, 1996.


\bibitem{Moklyachuk:2000}
\newblock M. P. Moklyachuk,
\newblock \emph{Robust procedures in time series analysis},
\newblock Theory of Stochastic Processes, vol. 6, no. 3-4, pp. 127-147, 2000.

\bibitem{Moklyachuk:2001}
\newblock M. P. Moklyachuk,
\newblock \emph{Game theory and convex optimization methods in robust estimation problems},
\newblock Theory of Stochastic Processes, vol. 7, no. 1-2, pp. 253-264, 2001.

\bibitem{Moklyachuk:2008}
\newblock M. P. Moklyachuk,
\newblock \emph{Robust estimations of  functionals of stochastic processes.},
\newblock Kyiv University, Kyiv, 2008.

\bibitem{Moklyachuk:2008nonsm}
\newblock M. P. Moklyachuk,
\newblock \emph{Nonsmooth analysis and optimization},
\newblock Kyiv University, Kyiv, 2008.

\bibitem{Moklyachuk:2015}
\newblock M. P. Moklyachuk,
\newblock \emph{Minimax-robust estimation problems for stationary stochastic sequences},
\newblock Statistics, Optimization \& Information Computing, vol. 3, no. 4, pp. 348 - 419, 2015.


\bibitem{Moklyachuk:2016}
\newblock M. Moklyachuk, and I. Golichenko,
\newblock \emph{Periodically correlated processes estimates},
\newblock LAP LAMBERT Academic Publishing, 2016.


\bibitem{Moklyachuk:2012}
\newblock M. Moklyachuk, and O. Masyutka,
\newblock \emph{Minimax-robust estimation technique for stationary stochastic processes},
\newblock LAP LAMBERT Academic Publishing, 2012.

\bibitem{MoklyachukYadrenko:1979}
\newblock M. P. Moklyachuk, and M. I. Yadrenko,
\newblock \emph{Linear statistical problems for homogeneous isotropic random fields on a sphere. I},
\newblock Theory of Probability and Mathematical Statistics, vol. 18, pp. 115-124, 1979.


\bibitem{MoklyachukYadrenko:1980}
\newblock M. P. Moklyachuk, and M. I. Yadrenko,
\newblock \emph{Linear statistical problems for homogeneous isotropic random fields on a sphere. II},
\newblock Theory of Probability and Mathematical Statistics, vol. 19, pp.129-139, 1980.


\bibitem{Muller}
\newblock C. M\"uller,
\newblock \emph{Spherical harmonics},
\newblock Lecture Notes in Mathematics 17. Berlin-Heidelberg-New York: Springer-Verlag, 1966.



\bibitem{Napolitano}
\newblock A. Napolitano,
\newblock \emph{Cyclostationarity: New trends and applications},
\newblock Signal Processing, vol. 120, pp. 385--408, 2016.


\bibitem{North}
\newblock G. R. North, and R. F. Cahalan,
\newblock \emph{Predictability in a solvable stochastic climate model},
\newblock J. Atmospheric Sciences 38, 504-513,, 1981.

\bibitem{Okamoto}
\newblock T. Okamoto, and W. Hu,
\newblock \emph{Angular trispectra of cmb temperature and polarization},
\newblock Phys. Rev., Vol. D66, p. 063008, 2002.



\bibitem{Rockafellar}
\newblock R. T. Rockafellar,
\newblock \emph{Convex Analysis},
\newblock Princeton University Press, 1997.


\bibitem{Rozanov}
\newblock Yu. A. Rozanov,
\newblock \emph{Stationary stochastic processes},
\newblock San Francisco-Cambridge-London-Amsterdam: Holden-Day, 1967.

\bibitem{Serpedin}
\newblock E. Serpedin, F. Panduru, I. Sari, and G. B. Giannakis,
\newblock  \emph{Bibliography on cyclostationarity},
\newblock Signal Processing, vol. 85, pp. 2233--2303, 2005.

\bibitem{SubbaRao2006}
\newblock T. Subba Rao and G. Terdik,
\newblock  \emph{Multivariate non-linear regression with applications},
\newblock In: P. Bertail, P. Doukhan, and P. Soulier (eds),
Dependence in Probability and Statistics. Springer Verlag, New York,
pp. 431-470, 2006.

\bibitem{SubbaRao2012}
\newblock T. Subba Rao and G. Terdik,
\newblock  \emph{Statistical analysis of spatio-temporal models and their applications},
\newblock In:  C. R. Rao (ed), Handbook of Statistics, Vol. 30, Elsevier B.V., pp. 521--541,  2012.


\bibitem{Terdik2015}
\newblock G. Terdik,
\newblock \emph{ Angular spectra for non-Gaussian isotropic fields},
\newblock Brazilian Journal of Probability and Statistics,
vol. 29, no. 4, pp. 833--865, 2015.


\bibitem{Vastola}
\newblock  K. S. Vastola and H. V. Poor,
\newblock  \emph{ An analysis of the effects of spectral uncertainty on Wiener filtering},
\newblock  Automatica, vol. 28, pp. 289--293, 1983.


\bibitem{Wiener}
\newblock  N. Wiener,
\newblock \emph{Extrapolation, Interpolation and Smoothing of Stationary Time Series. With Engineering Applications},
\newblock  The M. I. T. Press, Massachusetts Institute of Technology, Cambridge, Mass., 1966.

\bibitem{Yadrenko}
\newblock M. I. Yadrenko,
\newblock \emph{Spectral theory of random fields},
\newblock Optimization Software Inc. Publications Division, New York, 1983.


\bibitem{Yaglom:1987a}
\newblock A. M. Yaglom,
\newblock \emph{Correlation theory of stationary and related random functions. Vol. 1: Basic results},
\newblock Springer Series in Statistics, Springer-Verlag, New York etc., 1987.

\bibitem{Yaglom:1987b}
\newblock A. M. Yaglom,
\newblock \emph{Correlation theory of stationary and related random functions. Vol. 2: Suplementary notes and references},
\newblock Springer Series in Statistics, Springer-Verlag, New York etc., 1987.



\end{thebibliography}
\end{document}